\documentclass[12pt]{article}

\usepackage{amssymb}
\usepackage{amsfonts}
\usepackage{graphicx}
\usepackage{caption}
\def\F{\mathbb{F}}
\def\R{\mathbb{R}}
\def\C{\mathbb{C}}
\def\H{\mathbb{H}}
\def\OO{\mathbb{O}}
\def\Z{\mathbb{Z}}
\def\om{\omega}
\newtheorem{thm}{Theorem}

\def\cqfd{\hfill $\Box$\break\null}
\def\proof{\noindent{\it Proof.\  }}

\begin{document}

\title{Octonion multiplication and Heawood's map}
\author{Bruno S\'evennec\thanks{UMPA ENS-Lyon, UMR 5669 CNRS, 46 All\'ee d'Italie, 69364 Lyon cedex 07, France.}
}

\maketitle

Almost any article or book dealing with Cayley-Graves algebra $\OO$ of octonions
(to be recalled shortly) has a picture like the following

\begin{figure}[h]
\centering
\includegraphics[height=5cm]{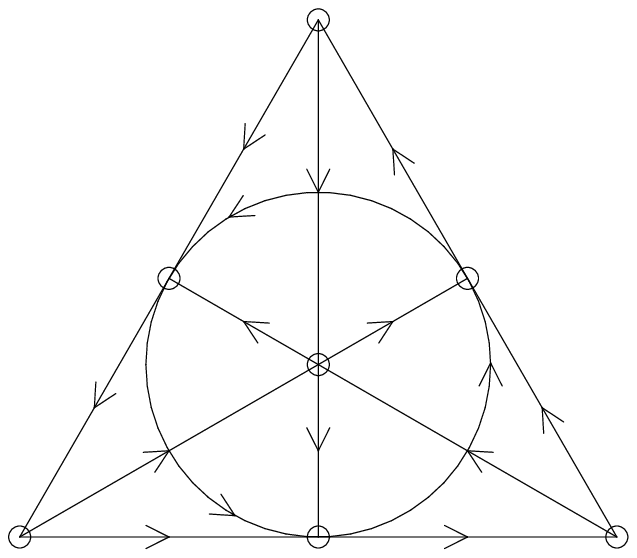}
\caption{}
\label{Fano}   
\end{figure}
\noindent
representing the so-called `Fano plane', which will be denoted
by $\Pi$, together with some cyclic ordering on each of its `lines'. It is a set of seven
points, in which seven three-point subsets called `lines' are specified, such that 
any two points are contained in a unique line, and any two lines
intersect in a unique point, giving a so-called (combinatorial) projective plane
\cite{vanLint-Wilson,Hughes-Piper}. In this cardinality, it is unique and may be seen as the
projective plane over the field $\F_2$. In the above figure, the points of $\Pi$ are represented
by vertices, midpoints of sides, and barycenter of an equilateral  triangle, and its `lines'
are either triples of aligned points (six sets) or the inscribed circle.
The arrows indicate a cyclic ordering on each line (six arrows are omitted for clarity).

The recipe to obtain the multiplication table of the seven imaginary basic
octonions $e_a$, $a\in\Pi$, from the above figure is then
\begin{itemize}
	\item For all $a$, $e_a^2=-1$.
	\item If $a,b,c$ are the cyclically ordered points of a line, then
	$$e_a e_b=e_c=-e_b e_a.$$
\end{itemize}
Finally one puts $$\OO=\R 1\oplus\bigoplus_{a\in\Pi}\R e_a \;\;.$$

Recall that $\OO$ is the largest $\R$-algebra with unit possessing a multiplicative
positive definite quadratic form $N:\OO\to\R$ 
and is in particular a division algebra.
The smaller ones are $\R,\C$, and $\H$ (Hamilton quaternion algebra),
and contrary to them $\OO$ is not associative, but only {\em alternative},
meaning for example that the associator $\{x,y,z\}=(xy)z-x(yz)$ is skew-symmetric
in $(x,y,z)$. Any three $e_a$'s with $a$'s on a line form together with $1$ a basis for an
isomorphic copy of $\H$ inside $\OO$.
One can find many interesting aspects of this algebra in \cite{Baez}
or \cite{Conway-Smith}, and also on the web at 
{\tt math.ucr.edu/home/baez/octonions/index.html}\ .

The purpose of this note is to present another graphical mnemonic for the
same multiplication table, which avoids the need to remember the seemingly 
arbitrary orientations on figure \ref{Fano}:

\begin{figure}[h]
\centering
\includegraphics[height=5cm]{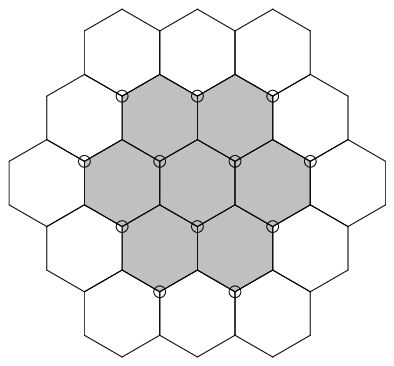}
\caption{}
\label{Flower} 
\end{figure}

This drawing is to be interpreted as follows. 
The union of the seven shaded hexagons is a fundamental domain
for a complex lattice $\Lambda\subset\C$, generated by `forward-right knight moves on an
hexagonal chessboard'. More mathematically, one can take $\Lambda=(2-\om)\Z[\om]$,
an ideal of norm (index)  $7$ in the ring $\Z[\om]$ of Eisenstein integers, $\om=(-1+\sqrt{-3})/2$.
Then the seven elements $0,\pm 1,\pm \om,\pm \om^2$ of $\Z[\om]$ form a complete set
of representatives of $\Z[\om]/\Lambda\simeq \Z/7\Z$, and the shaded hexagons
are their Dirichlet-Vorono\"{\i} cells relative to the lattice $\Z[\om]\subset\C$
(points which are closer to the considered point than to any other lattice points).

Going to the quotient torus $\C/\Lambda$, one sees on it seven hexagonal
regions, each of which is adjacent (along an edge) to the other six. Indeed this is obvious
for the one with center $0$, and translation invariance by $\Z[\om]/\Lambda$ implies
the claim. Dually, joining each hexagon's center $a\in\Z[\om]/\Lambda$ to
the centers $a\pm 1,a\pm \om,a\pm \om^2$ of its adjacent hexagons in $\C/\Lambda$, one obtains
an embedding of the complete graph $K_7$ on seven vertices in the torus.

In particular seven colors are necessary to color maps on a torus, and it is
not too difficult to show that they are also sufficient \cite[p. 434]{vanLint-Wilson}.
This last fact and the above map on the torus were discovered around 1890 by Heawood \cite{Heawood},
together with a gap in Kempe's 1879 proof of the four color theorem. 

Now any edge $\{a,b\}$ of this $K_7$ may be oriented
from $a$ to $b$ when $b-a\equiv 1,\om,\om^2 \;\mathrm{mod}\; \Lambda$ 
(equivalently if $b-a$ is a square mod $\Lambda$), and from $b$ to $a$
otherwise. As $1+\om+\om^2=0$, this defines a {\em cyclic orientation} 
on each set of three hexagons around a common vertex. The
circled vertices in figure \ref{Flower} are those where this orientation 
coincides with the complex (trigonometric) one (equivalently, those vertices where
the edges form a "Y").

The recipe to obtain the multiplication table of $\OO$ is now : associate to each hexagon on the
torus, {\it i.e.} to each vertex of $K_7$, or also to each element $a$ of $\Z[\om]/\Lambda\simeq \Z/7\Z$,
an imaginary basic octonion $e_a$ and declare that $e_a^2=-1$ and $e_a e_b = e_c =-e_b e_a$
whenever $a,b,c$ are cyclically ordered hexagons around a circled vertex.

To verify that the algebra thus defined is indeed $\OO$ one can check that the complex oriented
triangles in $K_7$, {\it i.e.} the triples of hexagons around a circled vertex, define the lines of 
a projective plane structure on the set of seven hexagons, and that the cyclic orientations
of the lines are the same as in figure \ref{Fano}, with appropriate identifications.

Alternatively, one may observe that the ring isomorphism $\Z[\om]/\Lambda\simeq \Z/7\Z$ sends
$\om\;\mathrm{mod}\; \Lambda$ to $2\;\mathrm{mod}\; 7$ and thus the multiplication rules
are the same as in \cite[p. 65]{Conway-Smith} or \cite[p. 150]{Baez},
namely --- with indices modulo 7 :

$$
\begin{array}{ccccc}
e_{n+1} e_{n+2} & = & e_{n+4} & = & -e_{n+2} e_{n+1} \\
e_{n+2} e_{n+4} & = & e_{n+1} & = & -e_{n+4} e_{n+2} \\
e_{n+4} e_{n+1} & = & e_{n+2} & = & -e_{n+1} e_{n+4}.
\end{array}
$$

\begin{figure}[h]
\centering
\includegraphics[height=5cm]{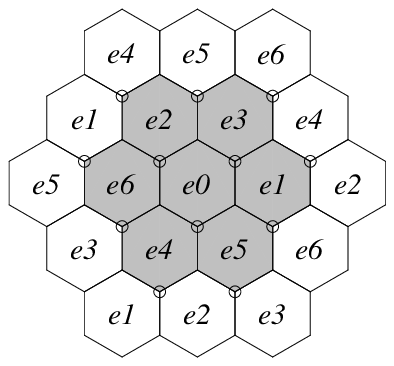}
\end{figure}

However, as this may seem a mere coincidence, we propose the following perhaps more conceptual statement.
\begin{thm}
Up to isomorphism, there is a unique orientation of the
graph $K_7$ for which oriented circuits of length three
define a combinatorial triangulation
of a surface (without boundary). This surface is a torus.
Moreover, given an orientation of this torus,
the triangles whose boundary orientation coincides with the graph orientation
define a (Fano) projective plane structure on the set of vertices. The cyclic orderings of the
lines in this plane coincide with those of figure \ref{Fano}.
\end{thm}
\proof 
Let us first consider a combinatorial triangulation of a surface $S$ having $K_7$ as its $1$-skeleton.
It has 21 edges, each belonging to exactly two triangles, so there are $t$ triangles, with $3t=42$,
{\it i.e.} $t=14$. But this gives for the Euler characteristic $\chi(S)=7-21+14=0$, hence $S$ is a torus
or a Klein bottle. We will exclude the second possibility 
(the proof is identical to that in \cite[p. 157]{Bollobas}, and is given here for convenience).

Observe that around each vertex $x$ of $K_7$,
the triangles containing $x$ together with the `adjacency along an edge from $x$' relation, 
must form a single cycle of length 6, because $S$ is a surface. 
Label $0$ some vertex and $1,2,\dots,6$ the others,
in trigonometric order around $0$. 
Then consider the universal cover $\widetilde{S}$ of $S$ together with
the lifted triangulation and labeling of vertices by $0,1,2,\dots,6$. 

This gives a labeling
of the vertices of the equilateral triangulation of the plane, such that all six labels different from $a$
appear around each vertex labelled $a$, in a cycle that depends only on $a$ (up to orientation, to allow for
the Klein bottle possibility).

\begin{figure}[h]
\centering
\includegraphics[height=5cm]{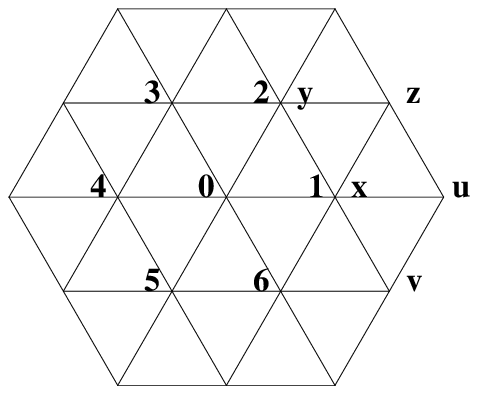}
\end{figure}
 
By the above choice, $1,2,\dots,6$ appear in cyclic order around
a $0$-labelled vertex (`first layer'); then consider the labelling possibilities for the vertices at
combinatorial distance $2$ from this $0$-vertex (`second layer') --- there are $12$ such vertices. It is easy
to see for instance that the second layer vertex $z$ adjacent to the first layer vertices $x,y$ labelled 
respectively $1$ and $2$ can only be labelled $4$ or $5$, since labels $0,1,2,3,6$ are excluded 
--- they already appear around $x$ or $y$. 

Assume that $z$ is labelled $5$ (the other choice is symmetrical).
Denote by $u$ the second layer vertex $u$ adjacent to $x,z$ and by $v$
the second layer vertex adjacent to $x$ and $u$ (thus $z,u,v$ appear in clockwise order around $x$).
Then $u$ can be labelled $3$ or $4$. If it were $4$, then $v$ is labelled $3$, and around
label $4$ one would have sequences of labels $5-1-3$ (from vertex $u$), and also $5-0-3$ (from the
first layer vertex labeled $4$). This contradicts the fact that all labels must appear in a single cycle
around each vertex of $S$. Thus $u$ is labelled $3$, and  $v$
is labelled $4$. It is now fairly easy to see that the labelling of second layer vertices
is uniquely determined (2,3,1,2,6,1,5,6,4 in clockwise order after vertex $v$), and that the deck
transformations of $\widetilde{S}$ are translations by vectors in the lattice $\Lambda=(2-\om)\Z[\om]$
of `forward right knight moves'. 

\begin{figure}[h]
\centering
\includegraphics[height=5cm]{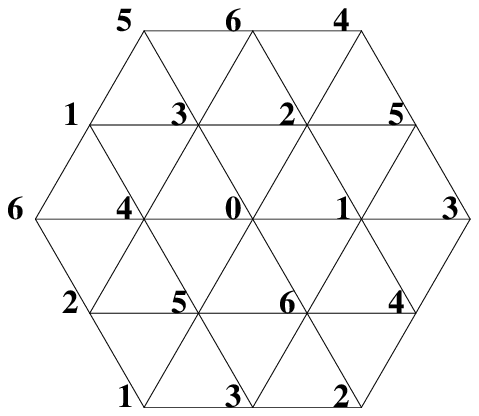}
\end{figure}

Note that in case $z$ is labelled $4$ instead, one obtains the mirror image
lattice $\overline{\Lambda}=(2-\om^2)\Z[\om]$ of `forward {\it left} knight moves'.

In particular $S$ is necessarily orientable (i.e. a torus), and up to (simplicial) isomorphism,
there is only one triangulation of a surface having $K_7$ as its $1$-skeleton.

To finish the proof of the first statement in the theorem, 
observe that the resulting $14$ triangles map is $2$-colorable
(equivalently, the dual graph is bipartite), for instance coloring in black triangles `pointing up'
in the universal covering plane. Orienting the edges of $K_7$ as boundaries of
black triangles (for some orientation of the torus) gives a solution to the first statement,
and the uniqueness (up to isomorphism) results from that of the triangulation and from consideration of
$-\mathrm{Id}$, which exchanges the two colors.

Checking that coherently oriented (alias black) triangles define the lines of
a projective plane is easy: they obviously form a set of  seven `lines', $3$-element subsets
of a $7$-element set of `points', such that any point is on exactly $3$ lines and any two points
are on a unique line. A simple counting argument shows that any two lines intersect then in a unique point.

To verify that the resulting cyclic orderings of lines correspond to those represented
in figure \ref{Fano}, one only has to see (by uniqueness)
that the latter verify the triangulation property
of the theorem. This is not difficult, taking symmetries into account.
Note that to obtain an oriented $K_7$, one must add six oriented edges to the
graph of figure \ref{Fano}, namely those ``closing''
the three sides and three medians. Then the reader will quickly check that any oriented
edge is contained in exactly two oriented circuits of length $3$, and that those oriented `triangles'
containing any vertex $x$ form a single cycle with respect to the relation of adjacency
along an edge from $x$.

\cqfd

\medskip

\noindent{\bf Remarks.\ }
i) The $K_7$ graph on the torus has still another link to the Fano projective plane $\Pi$.
Namely, its dual graph, union of the boundaries of the seven hexagons map, is the
incidence graph of $\Pi$: it is the bipartite graph with vertices the points and lines of $\Pi$,
with an edge for each membership. It is also known as the Heawood graph.
\medskip

\noindent
ii) Since it is easy to show that any embedding of $K_7$ in a closed surface 
of Euler characteristic $0$ has triangular faces, the above proof shows that $K_7$ embeds
uniquely (up to homeomorphism) in the torus, and also recovers the fact that 
it doesn't embed in the Klein bottle.

\end{document}